\documentclass{amsproc}%
\usepackage{graphicx}
\usepackage{amscd}
\usepackage{amsmath}
\usepackage{amsfonts}
\usepackage{amssymb}%
\theoremstyle{plain}

\numberwithin{equation}{section}

\begin{document}
\title[The zero-level centralizer in endomorphism algebras]{The zero-level centralizer in endomorphism algebras }
\author{Jen\H{o} Szigeti}
\address{Institute of Mathematics, University of Miskolc, Miskolc, Hungary 3515}
\email{jeno.szigeti@uni-miskolc.hu}
\author{Leon van Wyk}
\address{Department of Mathematical Sciences, Stellenbosch University\\
P/Bag X1, Matieland 7602, Stellenbosch, South Africa }
\email{LvW@sun.ac.za}
\thanks{The authors were supported by the National Research Foundation of South Africa
under Grant No.~UID 72375. Any opinion, findings and conclusions or
recommendations expressed in this material are those of the authors and
therefore the National Research Foundation does not accept any liability in
regard thereto.}
\subjclass[2010]{15A30, 15A27, 16D60, 16S50, 16U70}
\keywords{zero-level centralizer of a module endomorphism, nilpotent Jordan normal base}

\begin{abstract}
For an endomorphism $\varphi\in\mathrm{End}_{R}(M)$ of a left $R$-module
$_{R}M$ we investigate the structure and the polynomial identities of the
zero-level centralizer $\mathrm{Cen}_{0}(\varphi)$\ and the factor
\textrm{Cen}$(\varphi)/\mathrm{Cen}_{0}(\varphi)$. A double zero-centralizer
theorem for $\mathrm{Cen}_{0}(\mathrm{Cen}_{0}(\varphi))$ is also formulated.

\end{abstract}
\maketitle

\noindent1. INTRODUCTION

\bigskip

If $S$ is a ring (or algebra), then the centralizer \textrm{Cen}$(s)=\{u\in
S\mid us=su\}$ of an element $s\in S$ is a subring (subalgebra) of $S$. We
have $\mathrm{Cen}(s)=%
{\textstyle\bigcup\nolimits_{c\in\mathrm{LCen}(s)}}
\mathrm{Cen}_{c}(s)$, where $\mathrm{Cen}_{c}(s)=\{u\in S\mid us=su=c\}$ is
called the $c$-level centralizer and $\mathrm{LCen}(s)=\{c\in S\mid
\mathrm{Cen}_{c}(s)\neq\varnothing\}$ is a subring of $\mathrm{Cen}(s)$. The
zero-level centralizer $\mathrm{Cen}_{0}(s)=\{u\in S\mid us=su=0\}$ (or the
two-sided annihilator) of $s$ is an ideal of \textrm{Cen}$(s)$ and
$u+\mathrm{Cen}_{0}(s)\longmapsto us$ is a natural $\mathrm{Cen}%
(s)/\mathrm{Cen}_{0}(s)\longrightarrow\mathrm{LCen}(s)$ isomorphism of the
additive Abelian groups.

The aim of this paper is to investigate the zero-level centralizer
$\mathrm{Cen}_{0}(\varphi)$\ and the factor \textrm{Cen}$(\varphi
)/\mathrm{Cen}_{0}(\varphi)$\ for an element $\varphi$ in the endomorphism
ring $\mathrm{End}_{R}(M)$ of a left $R$-module $_{R}M$. Our treatment follows
the lines of [DSzW] and is heavily based on the results in [Sz] and [DSzW].
Thus we restrict our attention to the case of a finitely generated semisimple
$_{R}M$. First we focus on a nilpotent $\varphi$ and then we shall see that
for a non-nilpotent $\varphi$ the study of $\mathrm{Cen}_{0}(\varphi)$\ can be
reduced to the nilpotent case.

The authors were not able to find related results in the literature, in spite
of the fact that the objects of our investigations arise very naturally.
Surprisingly, the dimension formula for the zero-level centralizer of a square
matrix has not yet appeared in linear algebra books (e.g. [Ga,P,SuTy,TuA]).

In Section 2 we consider a fixed nilpotent Jordan normal base of $_{R}M$ with
respect to a given nilpotent $\varphi\in\mathrm{End}_{R}(M)$ and present all
the necessary prerequisites from [Sz] and [DSzW].

Section 3 is entirely devoted to the nilpotent case. Theorem 3.3 gives a
complete characterization of $\mathrm{Cen}_{0}(\varphi)$ and \textrm{Cen}%
$(\varphi)/\mathrm{Cen}_{0}(\varphi)$. If the base ring is local, then a more
accurate description of these algebras can be found in Theorem 3.4. Using 3.4
and the identities of certain subalgebras of a full matrix algebra over $R/J$,
in Theorems 3.6 and 3.8 we exhibit explicit polynomial identities for
$\mathrm{Cen}_{0}(\varphi)$ and \textrm{Cen}$(\varphi)/\mathrm{Cen}%
_{0}(\varphi)$, respectively.

In Section 4 we deal with the non-nilpotent case, a complete description of
$\mathrm{Cen}_{0}(\varphi)$ (as a particular ideal of an algebra of certain
invariant endomorphisms) can be found in Theorem 4.1. If $A\in M_{n}(K)$ is an
$n\times n$ matrix over a field $K$, then the mentioned dimension formula
$\dim_{K}\mathrm{Cen}_{0}(A)=\left[  \dim_{K}(\ker(A))\right]  ^{2}$ is an
immediate corollary of 4.1. Theorems 4.4 and 4.5 deal with the containment
relation $\mathrm{Cen}_{0}(\varphi)\subseteq\mathrm{Cen}_{0}(\sigma)$, where
$\sigma\in\mathrm{End}_{R}(M)$ is an other endomorphism. Since this
containment is equivalent to $\sigma\in\mathrm{Cen}_{0}(\mathrm{Cen}%
_{0}(\varphi))$, 4.4 and 4.5 can be considered as double zero-centralizer theorems.

\bigskip

\noindent2. PREREQUISITES

\bigskip

\noindent In order to provide a self-contained treatment, we collect some
notations, definitions and statements from [Sz] and [DSzW]. Let $Z(R)$ and
$J=J(R)$ denote the centre and the Jacobson radical of a ring $R$ (with
identity). Let $(z^{k})\vartriangleleft R[z]$ denote the ideal generated by
$z^{k}$ in the ring $R[z]$ of polynomials of the commuting indeterminate $z$.

\noindent For an $R$-endomorphism $\varphi:M\longrightarrow M$ of a (unitary)
left $R$-module $_{R}M$ a subset%
\[
\{x_{\gamma,i}\mid\gamma\in\Gamma,1\leq i\leq k_{\gamma}\}\subseteq M
\]
is called a nilpotent Jordan normal base of $_{R}M$\ with respect to $\varphi$
if each $R$-submodule $Rx_{\gamma,i}\leq M$ is simple, $\underset{\gamma
\in\Gamma,1\leq i\leq k_{\gamma}}{\oplus}Rx_{\gamma,i}=M$ is a direct sum,
$\varphi(x_{\gamma,i})=x_{\gamma,i+1}$, $\varphi(x_{\gamma,k_{\gamma}})=0$ for
all\textit{ }$\gamma\in\Gamma$, $1\leq i\leq k_{\gamma}$, and the set
$\{k_{\gamma}\mid\gamma\in\Gamma\}$ of integers is bounded. Now $\Gamma$\ is
called the set of (Jordan-) blocks and the size of the block $\gamma\in\Gamma$
is the integer $k_{\gamma}\geq1$.

\bigskip

\noindent\textbf{2.1.Theorem.}\textit{ Let }$\varphi\in\mathrm{End}_{R}%
(M)$\textit{ be an }$R$\textit{-endomorphism of a left }$R$\textit{-module
}$_{R}M$\textit{. Then the following are equivalent.}

\noindent1. $_{R}M$\textit{ is a semisimple left }$R$\textit{-module and
}$\varphi$\textit{\ is nilpotent of index }$n$.

\noindent2. \textit{There exists a nilpotent Jordan normal base }%
$X=\{x_{\gamma,i}\mid\gamma\in\Gamma,1\leq i\leq k_{\gamma}\}$\textit{ of}

\noindent\textit{ \ \ }$_{R}M$\textit{\ with respect to }$\varphi$\textit{
such that }$n=\max\{k_{\gamma}\mid\gamma\in\Gamma\}$.

\bigskip

\noindent\textbf{2.2.Theorem.}\textit{ Let }$\varphi\in\mathrm{End}_{R}%
(M)$\textit{ be a nilpotent }$R$\textit{-endomorphism of a finitely generated
semisimple left }$R$\textit{-module }$_{R}M$\textit{. If }%
\[
\{x_{\gamma,i}\mid\gamma\in\Gamma,1\leq i\leq k_{\gamma}\}\text{ \textit{and
}}\{y_{\delta,j}\mid\delta\in\Delta,1\leq j\leq l_{\delta}\}
\]
\textit{are nilpotent Jordan normal bases of }$_{R}M$\textit{\ with respect to
}$\varphi$\textit{, then }$\Gamma$\textit{ is finite and there exists a
bijection }$\pi:\Gamma\longrightarrow\Delta$\textit{ such that }$k_{\gamma
}=l_{\pi(\gamma)}$\textit{ for all }$\gamma\in\Gamma$\textit{. Thus the sizes
of the blocks of a nilpotent Jordan normal base are unique up to a permutation
of the blocks. We also have }$\ker(\varphi)=\underset{\gamma\in\Gamma}{\oplus
}Rx_{\gamma,k_{\gamma}}$\textit{ and hence }$\dim_{R}(\ker(\varphi
))=\left\vert \Gamma\right\vert $\textit{.}

\bigskip

\noindent If $\varphi\in\mathrm{End}_{R}(M)$ is an arbitrary $R$-endomorphism
of the left $R$-module $_{R}M$, then for $u\in M$ and $f(z)=a_{1}%
+a_{2}z+\cdots+a_{n+1}z^{n}\in R[z]$ (unusual use of indices!) the
multiplication%
\[
f(z)\ast u=a_{1}u+a_{2}\varphi(u)+\cdots+a_{n+1}\varphi^{n}(u)
\]
defines a natural left $R[z]$-module structure on $M$. This left action of
$R[z]$ on $M$ extends the left action of $R$ on $_{R}M$. For any
$R$-endomorphism $\psi\in\mathrm{End}_{R}(M)$ with $\psi\circ\varphi
=\varphi\circ\psi$ we have $\psi(f(z)\ast u)=f(z)\ast\psi(u)$ and hence
$\psi:M\longrightarrow M$ is an $R[z]$-endomorphism of the left $R[z]$-module
$_{R[z]}M$. On the other hand, if $\psi:M\longrightarrow M$ is an
$R[z]$-endomorphism of $_{R[z]}M$, then $\psi(\varphi(u))=\psi(z\ast
u)=z\ast\psi(u)=\varphi(\psi(u))$ implies that $\psi\circ\varphi=\varphi
\circ\psi$. Now $\mathrm{Cen}(\varphi)=\{\psi\mid\psi\in\mathrm{End}_{R}(M)$
and $\psi\circ\varphi=\varphi\circ\psi\}$ is a $Z(R)$-subalgebra of
$\mathrm{End}_{R}(M)$ and the argument above gives that $\mathrm{Cen}%
(\varphi)=\mathrm{End}_{R[z]}(M)$.

\noindent Henceforth $_{R}M$ is semisimple and we consider a fixed nilpotent
Jordan normal base
\[
X=\{x_{\gamma,i}\mid\gamma\in\Gamma,1\leq i\leq k_{\gamma}\}\subseteq M
\]
with respect to a given nilpotent $\varphi\in\mathrm{End}_{R}(M)$ of index
$n=\max\{k_{\gamma}\mid\gamma\in\Gamma\}$.

\noindent The $\Gamma$-copower $\amalg_{\gamma\in\Gamma}R[z]$\ is an ideal of
the $\Gamma$-direct power ring $(R[z])^{\Gamma}$ comprising all elements
$\mathbf{f}=(f_{\gamma}(z))_{\gamma\in\Gamma}$ with a finite set $\{\gamma
\in\Gamma\mid f_{\gamma}(z)\neq0\}$\ of non-zero coordinates. The copower
(power) has a natural $(R[z],R[z])$-bimodule structure. For an element
$\mathbf{f}=(f_{\gamma}(z))_{\gamma\in\Gamma}$ with $f_{\gamma}(z)=a_{\gamma
,1}+a_{\gamma,2}z+\cdots+a_{\gamma,n_{\gamma}+1}z^{n_{\gamma}}$ the formula%
\[
\Phi(\mathbf{f})=\underset{\gamma\in\Gamma,1\leq i\leq k_{\gamma}}{%
{\displaystyle\sum}
}a_{\gamma,i}x_{\gamma,i}=\underset{\gamma\in\Gamma}{%
{\displaystyle\sum}
}\left(  \underset{1\leq i\leq k_{\gamma}}{%
{\displaystyle\sum}
}a_{\gamma,i}\varphi^{i-1}(x_{\gamma,1})\right)  =\underset{\gamma\in\Gamma}{%
{\displaystyle\sum}
}f_{\gamma}(z)\ast x_{\gamma,1}%
\]
defines a function $\Phi:\amalg_{\gamma\in\Gamma}R[z]\rightarrow M$.

\bigskip

\noindent\textbf{2.3.Lemma.}\textit{ The function }$\Phi$\textit{\ is a
surjective left }$R[z]$\textit{-homomorphism. We have }$\varphi(\Phi
(\mathbf{f}))=\Phi(z\mathbf{f})$\textit{ for all }$\mathbf{f}\in\amalg
_{\gamma\in\Gamma}R[z]$\textit{ and the kernel}%
\[
\underset{\gamma\in\Gamma}{%
{\textstyle\coprod}
}J[z]+(z^{k_{\gamma}})\subseteq\ker(\Phi)\vartriangleleft_{l}\underset
{\gamma\in\Gamma}{%
{\textstyle\prod}
}R[z]
\]
\textit{is a left ideal of the power (and hence of the copower) ring. If }%
$R$\textit{ is a local ring (}$R/J$\textit{ is a division ring), then }%
$\amalg_{\gamma\in\Gamma}(J[z]+(z^{k_{\gamma}}))=\ker(\Phi)$\textit{.}

\bigskip

\noindent From now onward we also require that $_{R}M$ be finitely generated,
$m=\dim_{R}(\ker(\varphi))$, $\Gamma=\{1,2,\ldots,m\}$ and we assume that
$k_{1}\geq k_{2}\geq\ldots\geq k_{m}\geq1$ for the block sizes. Now
$\amalg_{\gamma\in\Gamma}R[z]=(R[z])^{\Gamma}$ and an element $\mathbf{f}%
=(f_{\gamma}(z))_{\gamma\in\Gamma}$ of $(R[z])^{\Gamma}$ is a $1\times m$
matrix (row vector) over $R[z]$. For an $m\times m$ matrix $\mathbf{P}%
=[p_{\delta,\gamma}(z)]$ in $\mathrm{M}_{m}(R[z])$ the matrix product%
\[
\mathbf{fP=}\underset{\delta\in\Gamma}{%
{\displaystyle\sum}
}f_{\delta}(z)\mathbf{p}_{\delta}%
\]
of $\mathbf{f}$ and $\mathbf{P}$ is a $1\times m$ matrix over $(R[z])^{\Gamma
}$, where $\mathbf{p}_{\delta}=(p_{\delta,\gamma}(z))_{\gamma\in\Gamma}$ is
the $\delta$-th row vector of $\mathbf{P}$ and%
\[
(\mathbf{fP})_{\gamma}=\underset{\delta\in\Gamma}{%
{\displaystyle\sum}
}f_{\delta}(z)p_{\delta,\gamma}(z).
\]
Consider the following subsets of $\mathrm{M}_{m}(R[z])$.%
\[
\mathcal{M}(X)=\{\mathbf{P}\in\mathrm{M}_{m}(R[z])\mid\mathbf{fP}\in\ker
(\Phi)\text{ for all }\mathbf{f}\in\ker(\Phi)\},
\]%
\[
\mathcal{I}(X)=\{\mathbf{P}\in\mathrm{M}_{m}(R[z])\mid\mathbf{P}%
=[p_{\delta,\gamma}(z)]\text{ and }p_{\delta,\gamma}(z)\in J[z]+(z^{k_{\gamma
}})\text{ for all }\delta,\gamma\in\Gamma\}=
\]%
\[
=\left[
\begin{array}
[c]{cccc}%
J[z]+(z^{k_{1}}) & J[z]+(z^{k_{2}}) & \cdots & J[z]+(z^{k_{m}})\\
J[z]+(z^{k_{1}}) & J[z]+(z^{k_{2}}) & \cdots & J[z]+(z^{k_{m}})\\
\vdots & \vdots & \ddots & \vdots\\
J[z]+(z^{k_{1}}) & J[z]+(z^{k_{2}}) & \cdots & J[z]+(z^{k_{m}})
\end{array}
\right]  ,
\]%
\[
\mathcal{N}(X)\!=\!\{\mathbf{P}\!\in\!\mathrm{M}_{m}(R[z])\!\mid
\!\mathbf{P}\!=\![p_{\delta,\gamma}(z)]\text{ and }z^{k_{\delta}}%
p_{\delta,\gamma}(z)\!\in\!J[z]\!+\!(z^{k_{\gamma}})\text{ for all }%
\delta,\gamma\!\in\!\Gamma\}.
\]
Note that $\mathcal{I}(X)$ and $\mathcal{N}(X)$ are $(R[z],R[z])$%
-sub-bimodules of $\mathrm{M}_{m}(R[z])$ in a natural way. For $\delta
,\gamma\in\Gamma$ let $k_{\delta,\gamma}=k_{\gamma}-k_{\delta}$ when $1\leq
k_{\delta}<k_{\gamma}\leq n$ and $k_{\delta,\gamma}=0$ otherwise. It can be
verified that the condition $z^{k_{\delta}}p_{\delta,\gamma}(z)\in
J[z]+(z^{k_{\gamma}})$ in the definition of $\mathcal{N}(X)$ is equivalent to
$p_{\delta,\gamma}(z)\in J[z]+(z^{k_{\delta,\gamma}})$ and so%
\[
\mathcal{N}(X)=\left[
\begin{array}
[c]{ccccc}%
R[z] & R[z] & R[z] & \cdots & R[z]\\
J[z]+(z^{k_{1}-k_{2}}) & R[z] & R[z] & \cdots & R[z]\\
J[z]+(z^{k_{1}-k_{3}}) & J[z]+(z^{k_{2}-k_{3}}) & R[z] & \cdots & R[z]\\
\vdots & \vdots & \vdots & \ddots & \vdots\\
J[z]+(z^{k_{1}-k_{m}}) & J[z]+(z^{k_{2}-k_{m}}) & J[z]+(z^{k_{3}-k_{m}}) &
\cdots & R[z]
\end{array}
\right]  .
\]

\bigskip

\noindent\textbf{2.4.Lemma.}\textit{ }$\mathcal{I}(X)\vartriangleleft
_{l}\mathrm{M}_{m}(R[z])$\textit{ is a left ideal, }$\mathcal{N}%
(X)\subseteq\mathrm{M}_{m}(R[z])$\textit{ is a subring, }$\mathcal{I}%
(X)\vartriangleleft\mathcal{N}(X)$\textit{ is an ideal and }$\mathcal{M}%
(X)$\textit{\ is a }$Z(R)$\textit{-subalgebra of }$\mathrm{M}_{m}%
(R[z])$\textit{. The ideal }$z\mathrm{M}_{m}(R[z])\vartriangleleft
\mathrm{M}_{m}(R[z])$\textit{ is nilpotent modulo }$\mathcal{I}(X)$\textit{
with }$\left(  z\mathrm{M}_{m}(R[z])\right)  ^{n}\subseteq\mathcal{I}%
(X)$.\textit{ If }$R$\textit{ is a local ring, then }$\mathcal{N}%
(X)=\mathcal{M}(X)$\textit{.}

\bigskip

\noindent\textbf{2.5.Theorem.}\textit{ Let }$\varphi\in\mathrm{End}_{R}%
(M)$\textit{ be a nilpotent }$R$\textit{-endomorphism of a finitely generated
semisimple left }$R$\textit{-module }$_{R}M$. \textit{For }$\mathbf{P}%
\in\mathcal{M}(X)$\textit{ and }$\mathbf{f}=(f_{\gamma}(z))_{\gamma\in\Gamma}%
$\textit{ in }$(R[z])^{\Gamma}$\textit{ the formula}%
\[
\psi_{\mathbf{P}}(\Phi(\mathbf{f}))=\Phi(\mathbf{fP})
\]
\textit{properly defines an }$R$\textit{-endomorphism }$\psi_{\mathbf{P}%
}\!:\!M\!\rightarrow\!M$\textit{ of }$_{R}M$\textit{ such that }%
$\psi_{\mathbf{P}}\!\circ\!\varphi\!=\!\varphi\!\circ\!\psi_{\mathbf{P}}%
$\textit{\ and the assignment }$\Lambda(\mathbf{P})=\psi_{\mathbf{P}}$\textit{
gives an }$\mathcal{M}(X)^{\mathrm{op}}\longrightarrow\mathrm{Cen}(\varphi
)$\textit{ homomorphism of }$Z(R)$\textit{-algebras. If }$\psi\circ
\varphi=\varphi\circ\psi$\textit{ holds for some }$\psi\in\mathrm{End}_{R}%
(M)$,\textit{ then there exists an }$m\times m$\textit{ matrix }$\mathbf{P}%
\in\mathcal{M}(X)$\textit{ such that }$\psi(\Phi(\mathbf{f}))=\Phi
(\mathbf{fP})$\textit{ for all }$\mathbf{f}=(f_{\gamma}(z))_{\gamma\in\Gamma}%
$\textit{ in }$(R[z])^{\Gamma}$\textit{. Thus }$\Lambda:\mathcal{M}%
(X)^{\mathrm{op}}\longrightarrow\mathrm{Cen}(\varphi)$\textit{ is surjective.}

\bigskip

\noindent\textbf{2.6.Lemma.}\textit{ }$\mathcal{I}(X)\subseteq\ker(\Lambda
)$\textit{ (}$\Lambda$\textit{ is defined in Theorem 2.5). If }$R$\textit{ is
a local ring then }$\mathcal{I}(X)=\ker(\Lambda)$\textit{.}

\bigskip

\noindent3. THE ZERO-LEVEL\ CENTRALIZER OF\ A NILPOTENT\ ENDOMORPHISM

\bigskip

\noindent We keep all settings from Section 2 and define the subsets of
$\mathrm{M}_{m}(R[z])\mathcal{\ }$as follows:%
\[
\mathcal{M}_{0}(X)=\{\mathbf{P}\in\mathcal{M}(X)\mid z\mathbf{fP}\in\ker
(\Phi)\text{ for all }\mathbf{f}\in(R[z])^{\Gamma}\},
\]%
\[
\mathcal{N}_{0}(X)\!=\!\{\mathbf{P}\in\mathrm{M}_{m}(R[z])\!\mid
\!\mathbf{P}=[p_{\delta,\gamma}(z)]\text{ and }p_{\delta,\gamma}(z)\!\in
J[z]+(z^{k_{\gamma}-1})\text{ for all }\delta,\gamma\!\in\Gamma\}.
\]
Since $p_{\delta,\gamma}(z)\in J[z]+(z^{k_{\gamma}-1})$ and $zp_{\delta
,\gamma}(z)\in J[z]+(z^{k_{\gamma}})$ are equivalent, we have%
\[
\mathcal{N}_{0}(X)=\left[
\begin{array}
[c]{cccc}%
J[z]+(z^{k_{1}-1}) & J[z]+(z^{k_{2}-1}) & \cdots & J[z]+(z^{k_{m}-1})\\
J[z]+(z^{k_{1}-1}) & J[z]+(z^{k_{2}-1}) & \cdots & J[z]+(z^{k_{m}-1})\\
\vdots & \vdots & \ddots & \vdots\\
J[z]+(z^{k_{1}-1}) & J[z]+(z^{k_{2}-1}) & \cdots & J[z]+(z^{k_{m}-1})
\end{array}
\right]  .
\]

\bigskip

\noindent\textbf{3.1.Lemma.}\textit{ }$\mathcal{I}(X)\subseteq\mathcal{N}%
_{0}(X)$\textit{, }$\left(  z\mathrm{M}_{m}(R[z])\right)  ^{n-1}%
\subseteq\mathcal{N}_{0}(X)$\textit{, }$\mathcal{N}_{0}(X)\vartriangleleft
_{l}\mathrm{M}_{m}(R[z])$\textit{ is a left ideal\thinspace and\thinspace
}$\mathcal{N}_{0}(X)\vartriangleleft\mathcal{N}(X)$\textit{\thinspace
is\thinspace an\thinspace ideal.\thinspace If\thinspace}$R$\textit{\thinspace
is\thinspace a\thinspace local\thinspace ring,\thinspace then\thinspace
}$\mathcal{N}_{0}(X)\!=\!\mathcal{M}_{0}(X)$\textit{.}

\bigskip

\noindent\textbf{Proof.} The containment $\mathcal{I}(X)\subseteq
\mathcal{N}_{0}(X)$ obviously holds and $\left(  z\mathrm{M}_{m}(R[z])\right)
^{n-1}\subseteq\mathcal{N}_{0}(X)$ is a consequence of $(z^{n-1}%
)\subseteq(z^{k_{\gamma}-1})$. Since the $\gamma$-th column of the matrices in
$\mathcal{N}_{0}(X)$ comes from a (left) ideal $J[z]+(z^{k_{\gamma}-1})$\ of
$R[z]$, we can see that $\mathcal{N}_{0}(X)$ is a left ideal of $\mathrm{M}%
_{m}(R[z])$.

\noindent If $\mathbf{P}\in\mathcal{N}_{0}(X)$ and $\mathbf{Q}\in
\mathcal{N}(X)$, then we have $zp_{\delta,\tau}(z)\in J[z]+(z^{k_{\tau}})$ and
$q_{\tau,\gamma}(z)\in J[z]+(z^{k_{\tau,\gamma}})$. Since $k_{\tau}%
+k_{\tau,\gamma}\geq k_{\gamma}$,\ it follows that $zp_{\delta,\tau}%
(z)q_{\tau,\gamma}(z)\in J[z]+(z^{k_{\gamma}})$. Thus $\mathbf{PQ}%
\in\mathcal{N}_{0}(X)$ and $\mathcal{N}_{0}(X)$ is an ideal of $\mathcal{N}%
(X)$.

\noindent If $R$\ is a local ring, then Lemma 2.3 gives that $\ker
(\Phi)=\amalg_{\gamma\in\Gamma}(J[z]+(z^{k_{\gamma}}))$. Let $\mathbf{1}%
_{\delta}$ denote the vector with $1$ in its $\delta$-coordinate and zeros in
all other places. If $\mathbf{P}\in\mathcal{M}_{0}(X)$, then $z\mathbf{1}%
_{\delta}\mathbf{P}\in\ker(\Phi)$ implies that $zp_{\delta,\gamma}(z)\in
J[z]+(z^{k_{\gamma}})$, whence $\mathbf{P}\in\mathcal{N}_{0}(X)$ follows. If
$\mathbf{P}\in\mathcal{N}_{0}(X)$ and $\mathbf{f}=(f_{\gamma}(z))_{\gamma
\in\Gamma}$ is in $(R[z])^{\Gamma}$, then $zp_{\delta,\gamma}(z)\in
J[z]+(z^{k_{\gamma}})$ implies that $zf_{\delta}(z)p_{\delta,\gamma}(z)\in
J[z]+(z^{k_{\gamma}})$ for all $\delta\in\Gamma$. Thus $z\mathbf{fP}\in
\ker(\Phi)$ and $\mathbf{P}\in\mathcal{M}_{0}(X)$ follows. $\square$

\bigskip

\noindent\textbf{3.2.Lemma.}\textit{\thinspace}$\ker(\Lambda)\!\subseteq
\!\mathcal{M}_{0}(X)$\textit{\thinspace and\thinspace for\thinspace
}$\mathbf{P\!}\in\mathit{\,}\mathcal{M}(X)$\textit{\thinspace the\thinspace
containments\thinspace}$\mathbf{P\!}\in\!\mathcal{M}_{0}(X)$\textit{\thinspace
and}

\noindent$\Lambda(\!\mathbf{P}\!)\!\in\!\mathrm{Cen}_{0}(\!\varphi\!)$\textit{
are\/equivalent.\/The\thinspace preimage }$\mathcal{M}_{0}(\!X\!)\!=\!\Lambda
^{-1}(\mathrm{Cen}_{0}(\!\varphi\!))\!\vartriangleleft\!\mathcal{M}%
(\!X\!)$\textit{\/\textit{i}s\/an\/ideal.}

\bigskip

\noindent\textbf{Proof.} The proof is based on the use of Lemma 2.3 and
Theorem 2.5.

\noindent If $\mathbf{P}\in\ker(\Lambda)$, then $\Lambda(\mathbf{P}%
)=\psi_{\mathbf{P}}=0$ gives that $\Phi(\mathbf{fP})=\psi_{\mathbf{P}}%
(\Phi(\mathbf{f}))=0$ for all $\mathbf{f}\in(R[z])^{\Gamma}$. Since
$\Phi:\amalg_{\gamma\in\Gamma}R[z]\rightarrow M$ is a left $R[z]$%
-homomorphism, $\Phi(z\mathbf{fP})=z\ast\Phi(\mathbf{fP})=0$ implies that
$z\mathbf{fP}\in\ker(\Phi)$. In view of $\ker(\Lambda)\subseteq\mathcal{M}%
(X)$, we deduce that $\mathbf{P}\in\mathcal{M}_{0}(X)$.

\noindent If $\mathbf{P}\in\mathcal{M}_{0}(X)$, then $\Lambda(\mathbf{P}%
)=\psi_{\mathbf{P}}$ and $\varphi(\psi_{\mathbf{P}}(\Phi(\mathbf{f}%
)))=\varphi(\Phi(\mathbf{fP}))=\Phi(z\mathbf{fP})=0$ for all $\mathbf{f}%
\in(R[z])^{\Gamma}$. Thus $\psi_{\mathbf{P}}\circ\varphi=\varphi\circ
\psi_{\mathbf{P}}=0$ and hence $\psi_{\mathbf{P}}\in\mathrm{Cen}_{0}(\varphi)$.

\noindent If $\Lambda(\mathbf{P})=\psi_{\mathbf{P}}$ is in $\mathrm{Cen}%
_{0}(\varphi)$, then $\varphi\circ\psi_{\mathbf{P}}=0$ and%
\[
\Phi(z\mathbf{fP})=\varphi(\Phi(\mathbf{fP}))=\varphi(\psi_{\mathbf{P}}%
(\Phi(\mathbf{f})))=0
\]
for all $\mathbf{f}\in(R[z])^{\Gamma}$. It follows that $\mathbf{P}%
\in\mathcal{M}_{0}(X)$.

\noindent Obviously, the preimage of the ideal $\mathrm{Cen}_{0}%
(\varphi)\vartriangleleft\mathrm{Cen}(\varphi)$ is also an ideal. $\square$

\bigskip

\noindent\textbf{3.3.Theorem.}\textit{ Let }$\varphi\in\mathrm{End}_{R}%
(M)$\textit{ be a nilpotent }$R$\textit{-endomorphism of a finitely generated
semisimple left }$R$\textit{-module }$_{R}M$. \textit{The map }$\Lambda
:\mathcal{M}(X)^{\mathrm{op}}\longrightarrow\mathrm{Cen}(\varphi)$\textit{
induces the following }$Z(R)$-\textit{isomorphisms for the factor algebras:}%
\[
\mathcal{M}_{0}(X)^{\mathrm{op}}/\ker(\Lambda)\cong\mathrm{Cen}_{0}%
(\varphi)\text{\textit{ and }}\mathcal{M}(X)^{\mathrm{op}}/\mathcal{M}%
_{0}(X)\cong\mathrm{Cen}(\varphi)/\mathrm{Cen}_{0}(\varphi).
\]

\bigskip

\noindent\textbf{Proof.} We have $\ker(\Lambda\upharpoonright\mathcal{M}%
_{0}(X))=\ker(\Lambda)$ and $\mathcal{M}_{0}(X)=\Lambda^{-1}(\mathrm{Cen}%
_{0}(\varphi))$\ by Lemma 3.2. Thus Theorem 2.5 ensures that the restricted
map $\Lambda\upharpoonright\mathcal{M}_{0}(X)$ is a surjective $\mathcal{M}%
_{0}(X)^{\mathrm{op}}\longrightarrow\mathrm{Cen}_{0}(\varphi)$ homomorphism of
$Z(R)$-algebras, whence $\mathcal{M}_{0}(X)^{\mathrm{op}}/\ker(\Lambda
)\cong\mathrm{Cen}_{0}(\varphi)$ follows.

\noindent In view of Lemma 3.2, the\ assignment%
\[
\mathbf{P}+\mathcal{M}_{0}(X)\longmapsto\Lambda(\mathbf{P})+\mathrm{Cen}%
_{0}(\varphi)
\]
is well-defined and gives an injective $\mathcal{M}(X)^{\mathrm{op}%
}/\mathcal{M}_{0}(X)\rightarrow\mathrm{Cen}(\varphi)/\mathrm{Cen}_{0}%
(\varphi)$\ homomorphism of $Z(R)$-algebras. The surjectivity of this
homomorphism is a consequence of the surjectivity of $\Lambda$ (see Theorem
2.5). $\square$

\bigskip

\noindent\textbf{3.4.Theorem.}\textit{ Let }$\varphi\in\mathrm{End}_{R}%
(M)$\textit{ be a nilpotent }$R$\textit{-endomorphism of a finitely generated
semisimple left }$R$\textit{-module }$_{R}M$. \textit{If }$R$\textit{\ is a
local ring, then the zero-level centralizer }$\mathrm{Cen}_{0}(\varphi
)$\textit{\ of }$\varphi$\textit{\ is isomorphic to the opposite of the factor
}$\mathcal{N}_{0}(X)/\mathcal{I}(X)$\textit{ as a }$Z(R)$\textit{-algebra:}%
\[
\mathrm{Cen}_{0}(\varphi)\cong\left(  \mathcal{N}_{0}(X)/\mathcal{I}%
(X)\right)  ^{\mathrm{op}}=\mathcal{N}_{0}(X)^{\mathrm{op}}/\mathcal{I}(X).
\]
\textit{We also have an isomorphism}%
\[
\mathrm{Cen}(\varphi)/\mathrm{Cen}_{0}(\varphi)\cong\left(  \mathcal{N}%
(X)/\mathcal{N}_{0}(X)\right)  ^{\mathrm{op}}=\mathcal{N}(X)^{\mathrm{op}%
}/\mathcal{N}_{0}(X)
\]
\textit{of the factor }$Z(R)$\textit{-algebras.}

\bigskip

\noindent\textbf{Proof.} Directly follows from Lemmas 2.4, 2.6, 3.1 and
Theorem 3.3. $\square$

\bigskip

\noindent Define a left ideal of $\mathrm{M}_{m}(R/J)$ as follows:%
\[
\mathcal{W}(X)=\{W=[w_{\delta,\gamma}]\mid w_{\delta,\gamma}\in R/J\text{ and
}w_{\delta,\gamma}=0\text{ if }k_{\gamma}\geq2\}.
\]
The assumption $k_{1}\geq k_{2}\geq\cdots\geq k_{m}\geq1$ ensures that%
\[
\mathcal{W}(X)=\left[
\begin{array}
[c]{cccccc}%
0 & \cdots & 0 & R/J & \cdots & R/J\\
0 & \ddots & 0 & R/J & \cdots & R/J\\
\vdots & \vdots & \ddots & \vdots & \vdots & \vdots\\
\vdots & \vdots & 0 & R/J & \cdots & R/J\\
\vdots & \vdots & \vdots & \vdots & \ddots & \vdots\\
0 & \cdots & 0 & R/J & \cdots & R/J
\end{array}
\right]  .
\]

\bigskip

\noindent\textbf{3.5.Lemma.}\textit{ }$\left(  \mathcal{N}_{0}(X)\cap
z\mathrm{M}_{m}(R[z])\right)  +\mathcal{I}(X)\vartriangleleft\mathcal{N}%
_{0}(X)$\textit{ is an ideal and there is a natural ring isomorphism}%
\[
\mathcal{N}_{0}(X)/(\left(  \mathcal{N}_{0}(X)\cap z\mathrm{M}_{m}%
(R[z])\right)  +\mathcal{I}(X))\cong\mathcal{W}(X)
\]
\textit{which is an }$(R,R)$\textit{-bimodule isomorphism at the same time.}

\bigskip

\noindent\textbf{Proof.} If $\mathbf{P}=[p_{\delta,\gamma}(z)]$ is in
$\mathcal{N}_{0}(X)$ and $p_{\delta,\gamma}(z)$\ has constant term
$u_{\delta,\gamma}\in R$, then%
\[
p_{\delta,\gamma}(z)-u_{\delta,\gamma}\in(J[z]+(z^{k_{\gamma}-1}))\cap(zR[z])
\]
and $k_{\gamma}\geq2$ implies that $u_{\delta,\gamma}\in J$. Thus
$[u_{\delta,\gamma}]\in\mathrm{M}_{m}(R)\cap\mathcal{N}_{0}(X)$ and%
\[
\mathbf{P}+(\left(  \mathcal{N}_{0}(X)\cap z\mathrm{M}_{m}(R[z])\right)
+\mathcal{I}(X))=[u_{\delta,\gamma}]+(\left(  \mathcal{N}_{0}(X)\cap
z\mathrm{M}_{m}(R[z])\right)  +\mathcal{I}(X))
\]
holds in $\mathcal{N}_{0}(X)/(\left(  \mathcal{N}_{0}(X)\cap z\mathrm{M}%
_{m}(R[z])\right)  +\mathcal{I}(X))$. The assignment%
\[
\mathbf{P}+(\left(  \mathcal{N}_{0}(X)\cap z\mathrm{M}_{m}(R[z])\right)
+\mathcal{I}(X))\longmapsto\lbrack u_{\delta,\gamma}+J]
\]
is well-defined and gives an%
\[
\mathcal{N}_{0}(X)/(\left(  \mathcal{N}_{0}(X)\cap z\mathrm{M}_{m}%
(R[z])\right)  +\mathcal{I}(X))\longrightarrow\mathcal{W}(X)
\]
isomorphism. $\square$

\bigskip

\noindent\textbf{3.6.Theorem.}\textit{\thinspace Let\thinspace}$R$%
\textit{\thinspace be\thinspace a\thinspace local\thinspace ring\thinspace
and\thinspace}$\varphi\in\mathrm{End}_{R}(M)$\textit{\thinspace be\thinspace
a\thinspace nilpotent\thinspace}$R$\textit{-endomorphism}

\noindent\textit{of\thinspace a\thinspace finitely\thinspace
generated\thinspace semisimple\thinspace left\thinspace}$R$\textit{-module
}$_{R}M$.\textit{\thinspace If\thinspace}$f_{i}(x_{1},\!\ldots\!,x_{r}%
)\!\in\!Z(R)\langle x_{1},\!\ldots\!,x_{r}\rangle$\textit{, }$1\!\leq
\!i\!\leq\!n$\textit{ and }$f_{i}=0$\textit{ are polynomial identities of the
right ideal }$\mathcal{W}(X)$\textit{ of }$\mathrm{M}_{m}^{\mathrm{op}}%
(R/J)$\textit{, then }$f_{1}f_{2}\cdots f_{n}=0$\textit{ is an identity of
}$\mathrm{Cen}_{0}(\varphi)$\textit{.}

\bigskip

\noindent\textbf{Proof.} Theorem 3.4 ensures that $\mathrm{Cen}_{0}%
(\varphi)\cong\mathcal{N}_{0}(X)^{\mathrm{op}}/\mathcal{I}(X)$ as
$Z(R)$-algebras, hence%
\[
Q=\left(  \left(  \mathcal{N}_{0}(X)\cap z\mathrm{M}_{m}(R[z])\right)
+\mathcal{I}(X)\right)  /\mathcal{I}(X)\vartriangleleft\mathcal{N}%
_{0}(X)/\mathcal{I}(X)
\]
can be viewed as an ideal of $\mathrm{Cen}_{0}(\varphi)$. The use of Lemma 3.5
gives%
\[
\mathrm{Cen}_{0}(\varphi)/Q\!\cong\!(\mathcal{N}_{0}(X)^{\mathrm{op}%
}\!/\mathcal{I}(X))\diagup\!Q\!\cong\!\mathcal{N}_{0}(X)^{\mathrm{op}%
}\!\diagup\!\left(  \mathcal{N}_{0}(X)\!\cap\!z\mathrm{M}_{m}(R[z])\right)
\!+\!\mathcal{I}(X)\!\cong\!\mathcal{W}(X)^{\mathrm{op}}\!.
\]
It follows that $f_{i}=0$ is an identity of $\mathrm{Cen}_{0}(\varphi)/Q$.
Thus $f_{i}(v_{1},\ldots,v_{r})\in Q$\ for all $v_{1},\ldots,v_{r}%
\in\mathrm{Cen}_{0}(\varphi)$, and so%
\[
f_{1}(v_{1},\ldots,v_{r})f_{2}(v_{1},\ldots,v_{r})\cdots f_{n}(v_{1}%
,\ldots,v_{r})\in Q^{n}.
\]
Since $\left(  z\mathrm{M}_{m}(R[z])\right)  ^{n}\subseteq\mathcal{I}(X)$ (see
Lemma 2.4) implies that $Q^{n}=\{0\}$, the proof is complete. $\square$

\bigskip

\noindent The assumption $k_{1}\geq k_{2}\geq\ldots\geq k_{m}\geq1$ ensures
that%
\[
\mathcal{U}_{0}(X)=\{U\in\mathrm{M}_{m}(R/J)\mid U=[u_{\delta,\gamma}]\text{
and }u_{\delta,\gamma}=0\text{ if }1\leq k_{\delta}<k_{\gamma}\text{ or
}k_{\gamma}=1\}.
\]
is a block upper triangular subalgebra of $\mathrm{M}_{m}(R/J)$. If
$[u_{\delta,\gamma}]\in\mathcal{U}_{0}(X)$ and $u_{\delta,\gamma}\neq0$ for
some $\delta,\gamma\in\Gamma$, then $2\leq k_{\gamma}\leq k_{\delta}$. Results
about the polynomial identities of block upper triangular matrix algebras can
be found in [GiZ].

\bigskip

\noindent\textbf{3.7.Lemma.} \textit{There is a natural ring isomorphism}%
\[
\mathcal{N}(X)/(\left(  \mathcal{N}(X)\cap z\mathrm{M}_{m}(R[z])\right)
+\mathcal{N}_{0}(X))\cong\mathcal{U}_{0}(X)
\]
\textit{which is an }$(R,R)$\textit{-bimodule isomorphism at the same time.}

\bigskip

\noindent\textbf{Proof.} For a matrix $\mathbf{P}=[p_{\delta,\gamma}(z)]$ in
$\mathcal{N}(X)$ consider the assignment%
\[
\mathbf{P}+\left(  \left(  \mathcal{N}(X)\cap z\mathrm{M}_{m}(R[z])\right)
+\mathcal{N}_{0}(X)\right)  \longmapsto\lbrack u_{\delta,\gamma}+J],
\]
where $u_{\delta,\gamma}\in R$ is defined as follows: $u_{\delta,\gamma}=0$ if
$k_{\gamma}=1$ and $u_{\delta,\gamma}$ is the constant term of $p_{\delta
,\gamma}(z)$\ if $k_{\gamma}\geq2$. Clearly, $[u_{\delta,\gamma}]\in
\mathrm{M}_{m}(R)\cap\mathcal{N}(X)$ and%
\[
\mathbf{P}+(\left(  \mathcal{N}(X)\cap z\mathrm{M}_{m}(R[z])\right)
+\mathcal{N}_{0}(X))=[u_{\delta,\gamma}]+(\left(  \mathcal{N}(X)\cap
z\mathrm{M}_{m}(R[z])\right)  +\mathcal{N}_{0}(X)).
\]
In view of the definitions of $\mathcal{N}_{0}(X)$ and $\mathcal{U}_{0}(X)$,
the above equality ensures that our assignment is a well-defined%
\[
\mathcal{N}(X)/(\left(  \mathcal{N}(X)\cap z\mathrm{M}_{m}(R[z])\right)
+\mathcal{N}_{0}(X))\longrightarrow\mathcal{U}_{0}(X)
\]
map providing the required isomomorphism. $\square$

\bigskip

\noindent\textbf{3.8.Theorem.}\textit{\thinspace Let\thinspace}$R$%
\textit{\thinspace be\thinspace a\thinspace local\thinspace ring\thinspace
and\thinspace}$\varphi\in\mathrm{End}_{R}(M)$\textit{\thinspace be\thinspace
a\thinspace nilpotent\thinspace}$R$\textit{-endomorphism}

\noindent\textit{of\thinspace a\thinspace finitely\thinspace
generated\thinspace semisimple\thinspace left\thinspace}$R$\textit{-module
}$_{R}M$.\textit{\thinspace If\thinspace}$f_{i}(x_{1},\!\ldots\!,x_{r}%
)\!\in\!Z(R)\langle x_{1},\!\ldots\!,x_{r}\rangle$\textit{, }$1\!\leq
\!i\!\leq\!n\!-\!1$\textit{ and }$f_{i}=0$\textit{ are polynomial identities
of the }$Z(R)$\textit{-subalgebra }$\mathcal{U}_{0}(X)$\textit{ of
}$\mathrm{M}_{m}^{\mathrm{op}}(R/J)$\textit{, then }$f_{1}f_{2}\cdots
f_{n-1}=0$\textit{ is an identity of the factor }$\mathrm{Cen}(\varphi
)/\mathrm{Cen}_{0}(\varphi)$\textit{.}

\bigskip

\noindent\textbf{Proof.} Theorem 3.4 ensures that $\mathrm{Cen}(\varphi
)/\mathrm{Cen}_{0}(\varphi)\cong\mathcal{N}(X)^{\mathrm{op}}/\mathcal{N}%
_{0}(X)$ as $Z(R)$-algebras, hence%
\[
L=\left(  \left(  \mathcal{N}(X)\cap z\mathrm{M}_{m}(R[z])\right)
+\mathcal{N}_{0}(X)\right)  /\mathcal{N}_{0}(X)\vartriangleleft\mathcal{N}%
(X)/\mathcal{N}_{0}(X)
\]
can be viewed as an ideal of $\mathrm{Cen}(\varphi)/\mathrm{Cen}_{0}(\varphi
)$. The use of Lemma 3.7\ gives%
\[
(\mathrm{Cen}(\varphi)/\mathrm{Cen}_{0}(\varphi))\diagup L\cong(\mathcal{N}%
(X)^{\mathrm{op}}/\mathcal{N}_{0}(X))\diagup L\cong
\]%
\[
\cong\mathcal{N}(X)^{\mathrm{op}}/(\left(  \mathcal{N}(X)\cap z\mathrm{M}%
_{m}(R[z])\right)  +\mathcal{N}_{0}(X))\cong\mathcal{U}_{0}(X)^{\mathrm{op}}.
\]
It follows that $f_{i}=0$ is an identity of $(\mathrm{Cen}(\varphi
)/\mathrm{Cen}_{0}(\varphi))\diagup L$. Thus $f_{i}(v_{1},\ldots,v_{r})\in
L$\ for all $v_{1},\ldots,v_{r}\in\mathrm{Cen}(\varphi)/\mathrm{Cen}%
_{0}(\varphi)$, and so%
\[
f_{1}(v_{1},\ldots,v_{r})f_{2}(v_{1},\ldots,v_{r})\cdots f_{n-1}(v_{1}%
,\ldots,v_{r})\in L^{n-1}.
\]
Since $\left(  z\mathrm{M}_{m}(R[z])\right)  ^{n-1}\subseteq\mathcal{N}%
_{0}(X)$ (see Lemma 3.1) implies that $L^{n-1}=\{0\}$, the proof is complete.
$\square$

\newpage

\noindent4. THE ZERO-LEVEL\ CENTRALIZER OF\ AN\ ARBITRARY ENDOMORPHISM

\bigskip

\noindent\textbf{4.1.Theorem.}\textit{ Let }$\varphi\in\mathrm{End}_{R}%
(M)$\textit{ be an }$R$\textit{-endomorphism of a finitely generated
semisimple left }$R$\textit{-module }$_{R}M$. \textit{Then there exist }%
$R$\textit{-submodules }$W_{1}$\textit{, }$W_{2}$\textit{ and }$V$\textit{\ of
}$M$\textit{ such that }$W=W_{1}\oplus W_{2}$\textit{ and }$M=V\oplus
W$\textit{ are direct products, }$\ker(\varphi)\subseteq W$\textit{, }%
$\varphi(W)=W_{2}$\textit{, }$\varphi(V)=V$\textit{, }$\dim_{R}(W_{1}%
)=\dim_{R}(\ker(\varphi))$\textit{, }$(\varphi\upharpoonright W)\in
\mathrm{End}_{R}(W)$\textit{ is nilpotent and for the zero-level centralizer
of }$\varphi$\textit{ we have }$\mathrm{Cen}_{0}(\varphi)\cong$\textit{\ }%
$T$\textit{, where}%
\[
T=\{\theta\in\mathrm{End}_{R}(W)\mid\theta(W_{1})\subseteq\ker(\varphi)\text{
and }\theta(W_{2})=\{0\}\}=\mathrm{Cen}_{0}(\varphi\upharpoonright W)
\]
\textit{is a left ideal of}%
\[
\mathrm{End}_{R}^{\ast}(W)=\{\alpha\in\mathrm{End}_{R}(W)\mid\alpha
(\ker(\varphi))\subseteq\ker(\varphi)\}
\]
\textit{and a right ideal of}%
\[
\mathrm{End}_{R}^{\ast\ast}(W)=\{\alpha\in\mathrm{End}_{R}(W)\mid\alpha
(W_{1}+\ker(\varphi))\subseteq W_{1}+\ker(\varphi)\text{\thinspace
and\thinspace}\alpha(W_{2})\subseteq W_{2}\}.
\]

\bigskip

\noindent\textbf{Proof.} The Fitting Lemma ensures the existence of an integer
$t\geq1$ such that $\mathrm{im}(\varphi^{t})\oplus\ker(\varphi^{t})=M$ is a
direct sum, where the (left) $R$-submodules%
\[
V=\mathrm{im}(\varphi^{t})=\mathrm{im}(\varphi^{t+1})=\cdots\text{ and }%
W=\ker(\varphi^{t})=\ker(\varphi^{t+1})=\cdots
\]
of $_{R}M$\ are uniquely determined by $\varphi$. Clearly, $\varphi(V)=V$ and
$\varphi(W)\subseteq W$ and the restricted map $(\varphi\upharpoonright
W)\in\mathrm{End}_{R}(W)$ is nilpotent of index $q\geq1$, where $\ker
(\varphi^{q-1})\neq\ker(\varphi^{q})=W$. Since $_{R}W$ is also finitely
generated and semisimple, Theorem 2.1 provides a nilpotent Jordan normal base
$X=\{x_{\gamma,i}\mid\gamma\in\Gamma,1\leq i\leq k_{\gamma}\}$ of $_{R}W$ with
respect to $\varphi\upharpoonright W$ (we have $x_{\gamma,k_{\gamma}+1}=0$ and
$q=\max\{k_{\gamma}\mid\gamma\in\Gamma\}$). Now $W_{1}\oplus W_{2}=W$ is a
direct sum, where%
\[
W_{1}=\underset{\gamma\in\Gamma}{\oplus}Rx_{\gamma,1}\text{ and }%
W_{2}=\underset{\gamma\in\Gamma,1\leq i\leq k_{\gamma}}{%
{\textstyle\bigoplus}
}Rx_{\gamma,i+1}\text{ }.
\]
Now we have $\ker(\varphi)\subseteq\ker(\varphi^{t})=W$ and $\ker
(\varphi)=\ker(\varphi\upharpoonright W)=\underset{\gamma\in\Gamma}{\oplus
}Rx_{\gamma,k_{\gamma}}$ by Theorem 2.2. It follows that%
\[
\dim_{R}(W_{1})=\left\vert \Gamma\right\vert =\dim_{R}(\ker(\varphi)).
\]
The definition of the nilpotent Jordan normal base ensures that $\varphi
(W)=W_{2}$.

\noindent If $\theta\in T$, then%
\[
\theta(\ker(\varphi))\subseteq\theta(W_{1}\oplus W_{2})=\theta(W_{1}%
)+\theta(W_{2})\subseteq\ker(\varphi)
\]
implies that $T$ is a left ideal of $\mathrm{End}_{R}^{\ast}(W)$ and a right
ideal of $\mathrm{End}_{R}^{\ast\ast}(W)$. Clearly, $T=\mathrm{Cen}%
_{0}(\varphi\upharpoonright W)$ is a consequence of $\varphi(W)=W_{2}$ and the
fact that $\theta(W)\subseteq\ker(\varphi)$ for all $\theta\in T$.

\noindent If $\alpha\in\mathrm{Cen}_{0}(\varphi)$, then $\alpha\circ\varphi=0$
implies that $\alpha(V)=\{0\}$ and $\alpha(x_{\gamma,i+1})=\alpha
(\varphi(x_{\gamma,i}))=0$ for $1\leq i\leq k_{\gamma}-1$. We also have
$\varphi\circ\alpha=0$, whence $\varphi(\alpha(x_{\gamma,1}))=0$ and
$\alpha(x_{\gamma,1})\in\ker(\varphi)$ follow. Thus $\alpha(W_{2})=\{0\}$,
$\alpha(W_{1})\subseteq\ker(\varphi)$ and the assignment $\alpha
\longmapsto\alpha\upharpoonright W$ obviously defines a $\mathrm{Cen}%
_{0}(\varphi)\longrightarrow T$ ring homomorphism.

\noindent If $\alpha,\beta\in\mathrm{Cen}_{0}(\varphi)$ and $\alpha
\upharpoonright W=\beta\upharpoonright W$, then $\alpha(V)=\beta(V)=\{0\}$ and
$V\oplus W=M$ ensure that $\alpha=\beta$ proving the injectivity of the above map.

\noindent If $\theta\in T$ and $\pi_{W}:V\oplus W\longrightarrow W$ is the
natural projection, then $\theta\circ\pi_{W}\in\mathrm{Cen}_{0}(\varphi)$.
Indeed, $\varphi\circ\theta\circ\pi_{W}=0$ is a consequence of $\theta
(W)\subseteq\ker(\varphi)$ and $\theta\circ\pi_{W}\circ\varphi=0$ is a
consequence of $\varphi(W)=W_{2}$ and $\theta(W_{2})=\{0\}$. Hence the
surjectivity of our assignment follows from $\theta\circ\pi_{W}\upharpoonright
W=\theta$. $\square$

\bigskip

\noindent\textbf{4.2.Corollary.}\textit{ Let }$A\in\mathrm{M}_{n}(K)$\textit{
be an }$n\times n$\textit{ matrix over a field }$K$\textit{, then the }%
$K$\textit{-dimension of the zero-level centralizer of }$A$\textit{ in
}$\mathrm{M}_{n}(K)$\textit{ is}%
\[
\dim_{K}\mathrm{Cen}_{0}(A)=\left[  \dim_{K}(\ker(A))\right]  ^{2}.
\]
\noindent\textbf{Proof.} Now $A\in\mathrm{End}_{K}(K^{n})$ and Theorem 4.1
ensures that $\mathrm{Cen}_{0}(A)\cong$\textit{\ }$T$, where%
\[
T=\{\theta\in\mathrm{End}_{K}(W)\mid\theta(W_{1})\subseteq\ker(A)\text{ and
}\theta(W_{2})=\{0\}\}.
\]
Our claim follows from the observation that the elements of $T$ and
$\mathrm{Hom}_{K}(W_{1},\ker(A))$ can be naturally identified and $\dim
_{K}(W_{1})=\dim_{K}(\ker(A))$. $\square$

\bigskip

\noindent\textbf{Remark.} Theorem 4.1 shows that the determination of the
zero-level centralizer can be reduced to the nilpotent case. This reduction
depends on the use of the Fitting Lemma.

\bigskip

\noindent\textbf{4.3.Lemma.}\textit{ Let }$\varphi,\sigma\in\mathrm{End}%
_{R}(M)$\textit{ be }$R$\textit{-endomorphisms of a finitely generated
semisimple left }$R$\textit{-module }$_{R}M$\textit{. If }$\mathrm{Cen}%
_{0}(\varphi)\subseteq\mathrm{Cen}_{0}(\sigma)$\textit{, then }$\mathrm{\ker
}(\varphi)\subseteq\mathrm{\ker}(\sigma)$\textit{ and }$\mathrm{im}%
(\sigma)\subseteq\mathrm{im}(\varphi)$\textit{.}

\bigskip

\noindent\textbf{Proof.} We use the proof of Theorem 4.1. If $\gamma\in\Gamma$
and $\pi_{\gamma}\in\mathrm{End}_{R}(M)$ denotes the natural%
\[
M=V\oplus W=V\oplus\left(  \underset{\delta\in\Gamma,1\leq i\leq k_{\delta}}{%
{\textstyle\bigoplus}
}Rx_{\delta,i}\right)  \longrightarrow Rx_{\gamma,k_{\gamma}}%
\]
projection, then $\pi_{\gamma}\circ\varphi^{k_{\gamma}-1}\in\mathrm{Cen}%
_{0}(\varphi)$. It follows that $\pi_{\gamma}\circ\varphi^{k_{\gamma}-1}%
\in\mathrm{Cen}_{0}(\sigma)$, whence we obtain that $\pi_{\gamma}\circ
\varphi^{k_{\gamma}-1}\circ\sigma=\sigma\circ\pi_{\gamma}\circ\varphi
^{k_{\gamma}-1}=0$. Since $\sigma(x_{\gamma,k_{\gamma}})=\sigma(\pi_{\gamma
}(\varphi^{k_{\gamma}-1}(x_{\gamma,k_{1}})))=0$, we have $x_{\gamma,k_{\gamma
}}\in\mathrm{\ker}(\sigma)$ for all $\gamma\in\Gamma$. Thus%
\[
\ker(\varphi)=\ker(\varphi\upharpoonright W)=\underset{\gamma\in\Gamma}%
{\oplus}Rx_{\gamma,k_{\gamma}}\subseteq\mathrm{\ker}(\sigma).
\]
The containment $\mathrm{im}(\sigma)\subseteq\ker(\pi_{\gamma}\circ
\varphi^{k_{\gamma}-1})$ is a consequence of $\pi_{\gamma}\circ\varphi
^{k_{\gamma}-1}\circ\sigma=0$, whence we obtain that $\mathrm{im}%
(\sigma)\subseteq\cap_{\gamma\in\Gamma}\ker(\pi_{\gamma}\circ\varphi
^{k_{\gamma}-1})$. It is straightforward to see that $\ker(\pi_{\gamma}%
\circ\varphi^{k_{\gamma}-1})=V\oplus W(\gamma)$ and%
\[
\underset{\gamma\in\Gamma}{%
{\textstyle\bigcap}
}(V\oplus W(\gamma))=V\oplus W_{2}=\varphi(V)+\varphi(W)=\varphi(V\oplus
W)=\mathrm{im}(\varphi),
\]
where%
\[
W(\gamma)=\underset{\delta\in\Gamma,1\leq i\leq k_{\delta},(\delta
,i)\neq(\gamma,1)}{%
{\textstyle\bigoplus}
}Rx_{\delta,i}.\text{ }\square
\]

\bigskip

\noindent\textbf{4.4.Theorem.}\textit{ Let }$\varphi,\sigma\in\mathrm{End}%
_{R}(M)$\textit{ be }$R$\textit{-endomorphisms of a finitely generated
semisimple left }$R$\textit{-module }$_{R}M$\textit{, then the following are
equivalent:}

\noindent1\textit{. }$\mathrm{Cen}_{0}(\varphi)\subseteq\mathrm{Cen}%
_{0}(\sigma)$\textit{,}

\noindent2.\textit{ }$\mathrm{\ker}(\varphi)\subseteq\mathrm{\ker}(\sigma
)$\textit{ and }$\mathrm{im}(\sigma)\subseteq\mathrm{im}(\varphi)$\textit{.}

\bigskip

\noindent\textbf{Proof.} In view of Lemma 4.3, it is enough to prove
$(2)\Longrightarrow(1)$. For an endomorphism $\tau\in\mathrm{Cen}_{0}%
(\varphi)$ we have $\tau\circ\varphi=\varphi\circ\tau=0$, whence
$\mathrm{im}(\sigma)\subseteq\mathrm{im}(\varphi)\subseteq\ker(\tau)$ and
$\mathrm{im}(\tau)\subseteq\mathrm{\ker}(\varphi)\subseteq\mathrm{\ker}%
(\sigma)$ follow. Thus we obtain that $\tau\circ\sigma=\sigma\circ\tau=0$. In
consequence we have $\tau\in\mathrm{Cen}_{0}(\sigma)$ and $\mathrm{Cen}%
_{0}(\varphi)\subseteq\mathrm{Cen}_{0}(\sigma)$ follows. $\square$

\bigskip

\noindent For a matrix $A\in\mathrm{M}_{n}(K)$ let $A^{\top}$ denote the
transpose of $A$.

\newpage

\noindent\textbf{4.5.Theorem.}\textit{ If }$A,B\in\mathrm{M}_{n}(K)$\textit{
are }$n\times n$ \textit{matrices over a field }$K$\textit{, then the
following are equivalent:}

\noindent1\textit{. }$\mathrm{Cen}_{0}(A)\subseteq\mathrm{Cen}_{0}%
(B)$\textit{,}

\noindent2.\textit{ }$\mathrm{\ker}(A)\subseteq\mathrm{\ker}(B)$\textit{ and
$\mathrm{\ker}$}$(A^{\top})\subseteq\ker(B^{\top})$\textit{,}

\noindent3.\textit{ }$\mathrm{im}(B)\subseteq\mathrm{im}(A)$\textit{ and
}$\mathrm{im}(B^{\top})\subseteq\mathrm{im}(A^{\top})$\textit{.}

\bigskip

\noindent\textbf{Proof.} $(1)\Longrightarrow(2)\&(3)$: For a matrix
$C\in\mathrm{Cen}_{0}(A^{\top})$ we have $CA^{\top}=A^{\top}C=0$ and $C^{\top
}\in\mathrm{Cen}_{0}(A)$ is a consequence of%
\[
AC^{\top}=(A^{\top})^{\top}C^{\top}=(CA^{\top})^{\top}=0=(A^{\top}C)^{\top
}=C^{\top}(A^{\top})^{\top}=C^{\top}A.
\]
Thus $C^{\top}\in\mathrm{Cen}_{0}(B)$ and a similar argument gives that
$C=(C^{\top})^{\top}\in\mathrm{Cen}_{0}(B^{\top})$. It follows that
$\mathrm{Cen}_{0}(A^{\top})\subseteq\mathrm{Cen}_{0}(B^{\top})$. The
application of Lemma 4.3 for the matrices $A,B,A^{\top},B^{\top}%
\in\mathrm{End}_{K}(K^{n})$ gives $\mathrm{\ker}(A)\subseteq\mathrm{\ker}(B)$,
$\mathrm{im}(B)\subseteq\mathrm{im}(A)$ and $\mathrm{\ker}(A^{\top}%
)\subseteq\mathrm{\ker}(B^{\top})$, $\mathrm{im}(B^{\top})\subseteq
\mathrm{im}(A^{\top})$.

\noindent$(2)\Longrightarrow(1)$: For a matrix $C\in\mathrm{Cen}_{0}(A)$ the
containment $\mathrm{im}(C)\subseteq\mathrm{\ker}(A)$ is a consequence of
$AC=0$ and $\mathrm{im}(C^{\top})\subseteq\mathrm{\ker}(A^{\top})$ is a
consequence of $A^{\top}C^{\top}=(CA)^{\top}=0$. Now $\mathrm{im}%
(C)\subseteq\mathrm{\ker}(B)$ implies that $BC=0$ and $\mathrm{im}(C^{\top
})\subseteq\mathrm{\ker}(B^{\top})$ implies that $CB=(B^{\top}C^{\top})^{\top
}=0$. Thus $C\in\mathrm{Cen}_{0}(B)$ and $\mathrm{Cen}_{0}(A)\subseteq
\mathrm{Cen}_{0}(B)$ follows.

\noindent$(3)\Longrightarrow(1)$: For a matrix $C\in\mathrm{Cen}_{0}(A)$ the
containment $\mathrm{im}(A)\subseteq\mathrm{\ker}(C)$ is a consequence of
$CA=0$ and $\mathrm{im}(A^{\top})\subseteq\mathrm{\ker}(C^{\top})$ is a
consequence of $C^{\top}A^{\top}=(AC)^{\top}=0$. Now $\mathrm{im}%
(B)\subseteq\mathrm{\ker}(C)$ implies that $CB=0$ and $\mathrm{im}(B^{\top
})\subseteq\mathrm{\ker}(C^{\top})$ implies that $BC=(C^{\top}B^{\top})^{\top
}=0$. Thus $C\in\mathrm{Cen}_{0}(B)$ and $\mathrm{Cen}_{0}(A)\subseteq
\mathrm{Cen}_{0}(B)$ follows. $\square$

\bigskip

\noindent ACKNOWLEDGMENT: The authors wish to thank P.N. Anh and L. Marki for
fruitful consultations.

\bigskip

\noindent REFERENCES

\bigskip

\noindent\lbrack DSzW] Drensky, V., Szigeti, J. and van Wyk, L.
\textit{Centralizers in endomorphism rings,} J. Algebra \textbf{324} (2010), 3378-3387.

\noindent\lbrack Ga] Gantmacher, F.R. \textit{The Theory of Matrices}, Chelsea
Publishing Co., New York, 2000.

\noindent\lbrack GiZ] Giambruno, A. and Zaicev, M. \textit{Polynomial
Identities and Asymptotic Methods,} Mathematical Surveys and Monographs
\textbf{122}, Amer.~Math.~Soc., Providence, Rhode Island, 2005.

\noindent\lbrack P] Prasolov, V.V. \textit{Problems and Theorems in Linear
Algebra}, Vol. \textbf{134} of Translation of Mathematical Monographs,
Amer.~Math.~Soc., Providence, Rhode Island, 1994.

\noindent\lbrack SuTy] Suprunenko, D.A. and Tyshkevich, R.I.\textit{
Commutative Matrices}, Academic Press, New York and London, 1968.

\noindent\lbrack Sz] Szigeti, J. \textit{Linear algebra in lattices and
nilpotent endomorphisms of semisimple modules}, J. Algebra \textbf{319}
(2008), 296--308.

\noindent\lbrack TuA] Turnbull, H.W. and Aitken, A.C. \textit{An Introduction
to the Theory of Canonical Matrices}, Dover Publications, 2004.

\end{document}